\documentclass[12pt]{amsart}

\usepackage{amsmath}
\usepackage{amscd}
\usepackage{amssymb}
\usepackage{latexsym}
\usepackage{stmaryrd} 
\DeclareMathAlphabet{\mathbbb}{U}{bbold}{m}{n}

\newtheorem{theorem}{Theorem}[section]
\newtheorem*{theorem*}{Theorem}
\newtheorem{lemma}[theorem]{Lemma}
\newtheorem*{lemma*}{Lemma}
\newtheorem{corollary}[theorem]{Corollary}
\newtheorem{proposition}[theorem]{Proposition}

\newtheorem{remark}[theorem]{Remark}
\newtheorem{definition}[theorem]{Definition}

\makeatletter
\def\revddots{\mathinner{\mkern1mu\raise\p@
\vbox{\kern7\p@\hbox{.}}\mkern2mu
\raise4\p@\hbox{.}\mkern2mu\raise7\p@\hbox{.}\mkern1mu}}
\makeatother 
\newcommand{\bgl}{\begin{equation}} 
\newcommand{\egl}{\end{equation}}
\newcommand{\bgloz}{\begin{equation*}} 
\newcommand{\egloz}{\end{equation*}}
\newcommand{\bgln}{\begin{eqnarray}} 
\newcommand{\egln}{\end{eqnarray}}
\newcommand{\bglnoz}{\begin{eqnarray*}} 
\newcommand{\eglnoz}{\end{eqnarray*}}
\newcommand{\btheo}{\begin{theorem}}
\newcommand{\etheo}{\end{theorem}}
\newcommand{\btheooz}{\begin{theorem*}}
\newcommand{\etheooz}{\end{theorem*}}
\newcommand{\blemma}{\begin{lemma}}
\newcommand{\elemma}{\end{lemma}}
\newcommand{\blemmaoz}{\begin{lemma*}}
\newcommand{\elemmaoz}{\end{lemma*}}
\newcommand{\bproof}{\begin{proof}}
\newcommand{\eproof}{\end{proof}}
\newcommand{\bbew}{\begin{beweis}}
\newcommand{\ebew}{\end{beweis}}
\newcommand{\bremark}{\begin{remark}\em}
\newcommand{\eremark}{\end{remark}}
\newcommand{\bdefin}{\begin{definition}}
\newcommand{\edefin}{\end{definition}}
\newcommand{\bprop}{\begin{proposition}}
\newcommand{\eprop}{\end{proposition}}
\newcommand{\bcor}{\begin{corollary}}
\newcommand{\ecor}{\end{corollary}}
\newcommand{\bfa}{\begin{cases}} 
\newcommand{\efa}{\end{cases}}
\newcommand{\cK}{\mathcal K}
\newcommand{\cL}{\mathcal L}

\newcommand{\cO}{\mathcal O}

\newcommand{\cT}{\mathcal T}

\def\Az{\mathbb{A}}
\def\Cz{\mathbb{C}}
\def\Fz{\mathbb{F}}

\def\Nz{\mathbb{N}}

\def\Rz{\mathbb{R}}
\def\Tz{\mathbb{T}}
\def\Zz{\mathbb{Z}}

\def\1z{\mathbbb{1}}
\newcommand{\fA}{\mathfrak A}

\newcommand{\an}[1]{``#1''} 
\newcommand{\ti}{\tilde}

\newcommand{\ri}{\rightarrow}
\newcommand{\lori}{\longrightarrow}

\newcommand{\ma}{\mapsto} 
\newcommand{\loma}{\longmapsto} 
\newcommand\into{\hookrightarrow} 

\def\SEMI{\mbox{$\times\kern-2pt\vrule height5pt width.6pt \kern3pt $}}

\newcommand{\id}{{\rm id}}

\newcommand{\de}{{\rm deg\,}} 

\newcommand{\img}{{\rm im\,}}
\renewcommand{\ker}{{\rm ker}\,}

\newcommand{\reg}{^\times} 
\newcommand{\pos}{_{>0}} 
\newcommand{\abs}[1]{\lvert#1\rvert} 
\newcommand{\defeq}{\mathrel{:=}} 

\newcommand{\dop}{\text{: }} 
\newcommand{\fuer}{\text{ for }} 
\newcommand{\falls}{\text{ if }} 
\newcommand{\fa}{\text{ for all }} 
\newcommand{\ilim}{\varinjlim} 
\newcommand{\plim}{\varprojlim} 
\newcommand{\chf}[1]{\1z_{\left[#1\right]}} 
\newcommand{\extalg}{\Lambda^* \,} 
\newcommand{\rte}{\overset{e}{\rtimes}} 
\newcommand{\lge}{\left\{} 
\newcommand{\rge}{\right\}} 
\newcommand{\lru}{\left(} 
\newcommand{\rru}{\right)} 
\newcommand{\leck}{\left[} 
\newcommand{\reck}{\right]} 
\newcommand{\lsp}{\left\langle} 
\newcommand{\rsp}{\right\rangle} 
\newcommand{\rukl}[1]{\lru #1 \rru} 
\newcommand{\eckl}[1]{\leck #1 \reck} 
\newcommand{\gekl}[1]{\lge #1 \rge} 
\newcommand{\spkl}[1]{\lsp #1 \rsp} 
\newcommand{\menge}[2]{\gekl{ #1 \dop #2 }} 
\begin{document}

\title{K-theory for ring C*-algebras attached to function fields}

\author{Joachim Cuntz and Xin Li}


\thanks{\scriptsize{Research supported by the Deutsche Forschungsgemeinschaft (SFB 478).}}

\thanks{\scriptsize{The second named author is supported by the Deutsche Telekom Stiftung.}}

\thanks{\scriptsize{This work was done in the context of the second named author's PhD project at the University of M\"unster}}

\begin{abstract}

We compute the K-theory of ring C*-algebras for polynomial rings over finite fields. The key ingredient is a duality theorem which we had obtained in a previous paper. It allows us to show that the K-theory of these algebras has a ring structure and to determine explicit generators. Our main result also reveals striking similarities between the number field case and the function field case.
 
\end{abstract}

\maketitle

\section{Introduction}

The theory of ring C*-algebras, initiated in \cite{Cun3}, has been developed in \cite{CuLi1}, \cite{Li} and \cite{CuLi2}. The present paper continues our work in \cite{CuLi2} where we studied ring C*-algebras associated to rings of integers in number fields. In \cite{CuLi2} we proved a duality theorem which was a key ingredient in the computation of the K-theory of these algebras. It allowed us to pass from the finite adele ring to the infinite one where we could use homotopy arguments to determine the K-theory.

In the present paper, we turn to the function field case. Our goal is to compute the K-theory of the ring C*-algebra for $\Fz_q[T]$ where $q$ is a prime power, i.e. $q = p^n$ for some prime number $p$. Since our duality theorem holds for arbitrary global fields (see \cite{CuLi2}), we can apply it to function fields as well. However, since in that case the infinite adele space is totally disconnected, we cannot hope for homotopy arguments.

Nevertheless, the duality theorem and the passage from the finite to the infinite adele space give us, in a somewhat unexpected way, a different handle on the computation of K-theory. It allows us to find explicit generators for the K-theory which have sufficiently nice properties. These generators are not visible in the representation over the finite adele space. At the same time, this explicit description reveals a ring structure on the K-theory. 

As the final result, we obtain that the K-theory for the ring C*-algebra of $\Fz_q[T]$ can be described as the tensor product over $\Zz$ of $\ti{K}_0(C^*(\Fz_q \reg))$ and the exterior $\Zz$-algebra over the torsion-free part of the multiplicative group $\Fz_q(T) \reg$, where $\ti{K}_0(C^*(\Fz_q \reg))$ is the reduced K-theory of $C^*(\Fz_q \reg)$ (i.e. the cokernel of the canonical map $K_0(\Cz) \ri K_0(C^*(\Fz_q \reg))$). This formula is compatible with the ring structure.

We proceed as follows: First of all, we recall the notion of ring C*-algebras. We also summarize the results of \cite{CuLi2} (Section \ref{review}). Then, we determine the K-theory for the ring C*-algebra of $\Fz_q[T]$: First, using the duality theorem, we reduce our problem to computing $K_*(C_0(\Fz_q((T))) \rtimes \Fz_q(T) \rtimes \Fz_q(T) \reg)$ (Section \ref{dual}). Secondly, we start with computing $K_*(C_0(\Fz_q((T))) \rtimes \Fz_q(T) \rtimes (\Fz_q \reg \times \spkl{T}))$. It turns out that we can find explicit generators, projections and unitaries, for the K-groups (Section \ref{explK}). The crucial point is that these projections and unitaries commute with all the remaining unitaries one still has to adjoin in order to pass from $C_0(\Fz_q((T))) \rtimes \Fz_q(T) \rtimes (\Fz_q \reg \times \spkl{T})$ to $C_0(\Fz_q((T))) \rtimes \Fz_q(T) \rtimes \Fz_q(T) \reg$ (Section \ref{communit}). Finally, the computation is completed by comparing our situation with commutative tori of suitable dimensions (Section \ref{final}). 

\section{Review}
\label{review}

Let $K$ be a global field and $R$ the ring of integers in $K$. The ring C*-algebra $\fA[R]$ is defined as follows: Consider the Hilbert space $\ell^2(R)$ with canonical orthonormal basis $\menge{\xi_r}{r \in R}$. Define additive shifts $U^a$ via $U^a(\xi_r) = \xi_{a+r}$ and multiplicative shift operators $S_b$ by $S_b(\xi_r) = \xi_{br}$ (for $b \neq 0$). These unitaries and isometries generate a C*-subalgebra of $\cL(\ell^2(R))$, the ring C*-algebra $\fA[R]$. 
This concrete C*-algebra admits several alternative descriptions. For our purposes, the following one is important (see \cite{CuLi1}, Remark 3 and Section 5): 

\btheo
\label{AMAfin}
$\fA[R]$ is Morita equivalent to $C_0(\Az_f) \rtimes K \rtimes K \reg$. 
\etheo

The crossed product is taken with respect to the canonical action of $K \rtimes K \reg$ on the finite adele ring $\Az_f$ of $K$ via affine transformations. 

Moreover, we proved the following duality result (see \cite{CuLi2}, Theorem 4.1 and Corollary 4.2): 
\btheo
\label{duality}
For every global field $K$, the crossed products $C_0(\Az_{\infty}) \rtimes K \rtimes K \reg$ and $C_0(\Az_f) \rtimes K \rtimes K \reg$ are Morita equivalent. 
\etheo
Here $\Az_{\infty}$ is the infinite adele ring of $K$. The crossed products arise from the natural actions of $K \rtimes K \reg$ on $\Az_{\infty}$ and $\Az_f$ via affine transformations. 

This duality theorem allowed us to use homotopy arguments to determine $K_*(\fA[R])$ for the ring of integers $R$ in a number field $K$. Our final result is (see \cite{CuLi2}, Section 6): 

Let $K$ be a number field with roots of unity $\mu = \gekl{\pm 1}$ and ring of integers $R$. Let $\# \gekl{v_{\Rz}}$ be the number of real places of $K$. There is a decomposition $K \reg = \mu \times \Gamma$, where $\Gamma$ is a free abelian group on infinitely many generators, such that the K-theory of the ring C*-algebra of $R$ can be described as follows: 
\btheo
\label{Knf}
\bgloz
  K_*(\fA[R]) \cong 
  \bfa 
    K_0(C^*(\mu)) \otimes_{\Zz} \extalg (\Gamma) \falls \# \gekl{v_{\Rz}} = 0 \\
    \extalg (\Gamma) \falls \# \gekl{v_{\Rz}} \text{ is odd} \\
    \extalg (\Gamma) \oplus ((\Zz / 2 \Zz) \otimes_{\Zz} \extalg (\Gamma)) \falls \# \gekl{v_{\Rz}} \text{ is even and at least } 2.
  \efa
\egloz
\etheo
This isomorphism is meant as an isomorphism between $\Zz / 2 \Zz$-graded abelian groups. Here $K_*(\fA[R])$ is the canonically graded group $K_0(\fA[R]) \oplus K_1(\fA[R])$, $K_0(C^*(\mu))$ and $\Zz / 2 \Zz$ are trivially graded, the exterior $\Zz$-algebra $\extalg (\Gamma)$ is canonically graded and we consider graded tensor products. 

\section{Applying the duality theorem}
\label{dual}

Now we turn to function fields. Let us consider the case $K=\Fz_q(T)$ and $R=\Fz_q[T]$ for a prime power $q$. Our goal is to determine the K-theory of $\fA[R]$. By Theorem \ref{AMAfin}, we know that $\fA[R] \sim_M C_0(\Az_f) \rtimes K \rtimes K \reg$. Moreover, Theorem \ref{duality} yields $C_0(\Az_{\infty}) \rtimes K \rtimes K \reg \sim_M C_0(\Az_f) \rtimes K \rtimes K \reg$. Thus we have to compute the K-theory of $C_0(\Az_{\infty}) \rtimes K \rtimes K \reg$. 

It is our convention that the infinite adele ring over $K = \Fz_q(T)$ is given by 
\bgloz
  \Az_{\infty} \cong \Fz_q ((T)) = \menge{\sum_{i=n}^{\infty} a_i T^i}{n \in \Zz, a_i \in \Fz_q}.
\egloz 
$\Fz_q ((T))$ is a locally compact field with respect to the valuation 
\bgloz
  \abs{\sum_{i=n}^{\infty} a_i T^i}= q^{-n} \falls a_n \neq 0.
\egloz
Moreover, to form the crossed product $C_0(\Az_{\infty}) \rtimes K \rtimes K \reg$, we also need to know how $K$ sits inside $\Az_{\infty}$. The embedding $K \into \Az_{\infty}$ is determined by 
\bgloz
  K \supseteq \Fz_q [T] \ni a(T) \ma a(T^{-1}) \in \Fz_q ((T)) \cong \Az_{\infty} 
\egloz
(it is our convention that the infinite place of $K$ is given by the valuation $\abs{a/b}_{\infty} = q^{\de(a)-\de(b)}$ for $a \in \Fz_q [T]$, $b \in \Fz_q [T] \reg$).

Let $\ti{v}^a$, $\ti{t}_b$ be the unitaries in the multiplier algebra of $C_0(\Az_{\infty}) \rtimes K \rtimes K \reg$ which implement the additive and the multiplicative action, respectively. In other words, we have 
\bgloz
  \ti{v}^a \ti{t}_b f \ti{t}_b^* (\ti{v}^a)^* = f(\sigma(b)^{-1}(\sqcup - \sigma(a)))
\egloz
for every $f \in C_0(\Fz_q ((T)))$, where $\sigma$ is the ring isomorphism 
\bgloz
  \Fz_q (T) \ri \Fz_q (T); a(T) \ma a(T^{-1}). 
\egloz 
We observe that we can equally well consider the crossed product associated to the canonical embedding $\Fz_q (T) \into \Fz_q ((T))$. We denote this crossed product by $C_0(\Fz_q ((T))) \rtimes \Fz_q (T) \rtimes \Fz_q (T) \reg$ and let $v^a$, $t_b$ be the canonical unitaries in the multiplier algebra of this crossed product corresponding to addition and multiplication, respectively. We can identify $C_0(\Fz_q ((T))) \rtimes \Fz_q (T) \rtimes \Fz_q (T) \reg$ and $C_0(\Az_{\infty}) \rtimes K \rtimes K \reg$ via $f v^a t_b \ma f \ti{v}^{\sigma(a)} \ti{t}_{\sigma(b)}$. To be more precise, this homomorphism identifies the *-algebras $C_c(\Fz_q (T) \rtimes \Fz_q (T) \reg,C_0(\Fz_q ((T))))$ and $C_c(K \rtimes K \reg, C_0(\Az_{\infty}))$ viewed as *-subalgebras of $C_0(\Fz_q ((T))) \rtimes \Fz_q (T) \rtimes \Fz_q (T) \reg$ or $C_0(\Az_{\infty}) \rtimes K \rtimes K \reg$, respectively. Furthermore, this map is isometric with respect to the $\ell^1$-norms, so that it extends to an isomorphism of the crossed products. 

Thus our task is to determine the K-theory of $C_0(\Fz_q ((T))) \rtimes \Fz_q (T) \rtimes \Fz_q (T) \reg$. 

\section{Notations}
\label{not}

In the following, let $\chf{X}$ be the characteristic function of a subset $X$ in $\Fz_q ((T))$. In particular, the ring of power series $\Fz_q [[T]] = \menge{\sum_{i=0}^{\infty} a_i T^i}{a_i \in \Fz_q}$ sits inside $\Fz_q ((T))$, and we denote by $\1z_n$ the characteristic function $\chf{T^n \cdot \Fz_q [[T]]}$. The characteristic function of $\Fz_q [[T]]$ is denoted by $\1z$ (i.e. $\1z \defeq \1z_0$). Since the subset $T^n \cdot \Fz_q [[T]]$ is closed and open in $\Fz_q ((T))$, the functions $\1z_n$ and $\1z$ lie in $C_0(\Fz_q ((T)))$. 

Moreover, as we already had above, let $v^a$, $t_b$ be the canonical unitaries in the multiplier algebra of $C_0(\Fz_q ((T))) \rtimes \Fz_q (T) \rtimes \Fz_q (T) \reg$ implementing the additive or the multiplicative action, respectively. 

Furthermore, let $f_1, f_2, f_3, \dotsc$ be an enumeration of the irreducible polynomials in $\Fz_q [T]$ with constant term $1$, i.e. $f_i \in 1 + T \cdot \Fz_q [T]$. Let  $\Gamma$ be the subgroup of $\Fz_q(T) \reg$ generated by the polynomials $T$ and $f_1, f_2, f_3, \dotsc$. $\Gamma$ is a free abelian group, and free generators are precisely given by $T, f_1, f_2, f_3, \dotsc$. We have the decomposition $\Fz_q(T) \reg = \Fz_q \reg \times \Gamma$. Let $\Gamma_m \defeq \spkl{T, f_1, \dotsc, f_m}$. 

We will determine $K_*(C_0(\Fz_q ((T))) \rtimes \Fz_q (T) \rtimes \Fz_q (T) \reg)$ step by step, so it will be helpful to choose appropriate C*-subalgebras. Let 
\bgloz
  A_{-1} \defeq C_0(\Fz_q ((T))) \rtimes \Fz_q (T) \rtimes \Fz_q \reg
\egloz 
and 
\bgloz
A_m \defeq C_0(\Fz_q ((T))) \rtimes \Fz_q (T) \rtimes (\Fz_q \reg \times \Gamma_m) \fa m \in \Zz_{\geq 0}.
\egloz
If $\mu$ denotes the multiplicative action of $\Gamma$ on $A_{-1}$, then we have $A_0 \cong A_{-1} \rtimes_{\mu_T} \Zz$ and $A_m \cong A_{m-1} \rtimes_{\mu_{f_m}} \Zz$ for all $m \in \Zz_{\geq 0}$. Finally, $C_0(\Fz_q ((T))) \rtimes \Fz_q (T) \rtimes \Fz_q (T) \reg$ is isomorphic to $\ilim A_m$ with respect to the canonical maps $A_{m-1} \ri A_m$. Thus we have to determine the K-theory of $A_m$ for each $m$. 

\section{Explicit generators for K-theory}
\label{explK}

The first step is to determine the K-theory of $A_{-1}$. It turns out that $A_{-1}$ is approximately finite dimensional, so we just have to find a suitable description of $A_{-1}$ as an inductive limit of finite dimensional C*-algebras to compute its K-theory. 

\subsection{Filtrations}

Let $\mu_T$ be the endomorphism of $C(\Fz_q [[T]]) \rtimes (\Gamma \cdot \Fz_q [T]) \rtimes \Fz_q \reg$ induced by multiplication with $T$. $\mu_T$ is given by 
\bgloz
  \mu_T(f v^a t_b) = \rukl{f(T^{-1} \sqcup) \cdot \1z(T^{-1} \sqcup)} v^{aT} t_b. 
\egloz

\blemma
\label{f1}
$A_{-1}$ can be identified with the inductive limit of the system 
\bgloz
  \dotso \overset{\mu_T}{\ri} C(\Fz_q [[T]]) \rtimes (\Gamma \cdot \Fz_q [T]) \rtimes \Fz_q \reg 
  \overset{\mu_T}{\ri} C(\Fz_q [[T]]) \rtimes (\Gamma \cdot \Fz_q [T]) \rtimes \Fz_q \reg \overset{\mu_T}{\ri} \dotso
\egloz
\elemma

\bproof
The idea is that going over to this inductive limit corresponds to formally inverting $\mu_T$. 

To prove the claim, consider for each $n \in \Zz_{\geq 0}$ the homomorphism 
\bglnoz
  C(\Fz_q [[T]]) \rtimes (\Gamma \cdot \Fz_q [T]) \rtimes \Fz_q \reg &\ri& A_{-1} \\
  f v^a t_b &\ma& f(T^n \sqcup) v^{a/T^n} t_b.
\eglnoz
This family of homomorphisms is compatible with $\mu_T$ and thus gives rise to a homomorphism 
\bgloz
\ilim \gekl{C(\Fz_q [[T]]) \rtimes (\Gamma \cdot \Fz_q [T]) \rtimes \Fz_q \reg; \mu_T} \ri A_{-1}.
\egloz
This homomorphism is clearly surjective. To see injectivity, consider for each $n \in \Zz_{\geq 0}$ the commutative square 
\bgloz
  \begin{CD}
  C(\Fz_q [[T]]) \rtimes (\Gamma \cdot \Fz_q [T]) \rtimes \Fz_q \reg @>>> A_{-1} \\
  @VVV @VVV \\
  C(\Fz_q [[T]]) @>>> C_0(\Fz_q ((T)))
  \end{CD}
\egloz
where the upper horizontal arrow is the homomorphism introduced above (for the $n$ we have chosen) and the lower horizontal arrow is given by $f \ma f(T^n \sqcup)$. The vertical arrows are the canonical faithful conditional expectations. They exist because we are dealing with discrete amenable groups. As the lower horizontal homomorphism is clearly injective, the upper one has to be so as well. This proves injectivity for each $n$ and thus for the induced homomorphism on the inductive limit. 
\eproof

Now, for every $d \in \Gamma$ let $\mu_d$ be the endomorphism of $C(\Fz_q [[T]]) \rtimes \Fz_q [T] \rtimes \Fz_q \reg$ induced by multiplication with $d$. It is given by $f v^a t_b \ma f(d^{-1} \sqcup) v^{da} t_b$. We have 

\blemma
\label{f2}
\bgloz
  C(\Fz_q [[T]]) \rtimes (\Gamma \cdot \Fz_q [T]) \rtimes \Fz_q \reg \cong \ilim_{d \in \Gamma} \gekl{C(\Fz_q [[T]]) \rtimes \Fz_q [T] \rtimes \Fz_q \reg; \mu_d}.
\egloz
\elemma

\bproof
This can be proven analogously to the previous lemma. 
\eproof

Moreover, let $(\Fz_q [T])^{(n)}$ be the additive subgroup $\menge{a_0 + \dotsb + a_n T^n}{a_i \in \Fz_q}$ of $\Fz_q [T]$. For each $n$ in $\Zz \pos$, we can identify $(\Fz_q [T])^{(n-1)}$ and $\Fz_q [T] / T^n \cdot \Fz_q [T]$ as additive groups. Thus $(\Fz_q [T])^{(n-1)}$ acts additively on $C(\Fz_q [T] / T^n \cdot \Fz_q [T])$. This additive action and the multiplicative action of $\Fz_q \reg$ give rise to the crossed product 
\bgloz
  C(\Fz_q [T] / T^n \cdot \Fz_q [T]) \rtimes (\Fz_q [T])^{(n-1)} \rtimes \Fz_q \reg. 
\egloz
Let $\iota_{n,n+1}$ be the homomorphism 
\bgloz
  C(\Fz_q [T] / T^n \cdot \Fz_q [T]) \rtimes (\Fz_q [T])^{(n-1)} \rtimes \Fz_q \reg \ri C(\Fz_q [T] / T^{n+1} \cdot \Fz_q [T]) \rtimes (\Fz_q [T])^{(n)} \rtimes \Fz_q \reg
\egloz
given by $g v^a t_b \ma (g \circ \pi_{n+1,n}) v^a t_b$ with the canonical projection $\pi_{n+1,n}$ from $\Fz_q [T] / T^{n+1} \cdot \Fz_q [T]$ onto $\Fz_q [T] / T^n \cdot \Fz_q [T]$. We have 
\blemma
\label{f3}
\bgloz
  C(\Fz_q [[T]]) \rtimes \Fz_q [T] \rtimes \Fz_q \reg \cong \ilim_n \gekl{C(\Fz_q [T] / T^n \cdot \Fz_q [T]) \rtimes (\Fz_q [T])^{(n-1)} \rtimes \Fz_q \reg;\iota_{n,n+1}}
\egloz
\elemma
\bproof
Again, the proof is analogous to the one of Lemma \ref{f1}. The point is that $\Fz_q [[T]]$ can be identified with 
\bgloz
  \plim_n \gekl{\Fz_q [T] / T^{n+1} \cdot \Fz_q [T]; \pi_{n+1,n}}
\egloz
both algebraically and topologically.
\eproof

\subsection{Explicit generators for the K-groups}

The preceding filtrations allow us to compute the K-theory of $A_{-1}$. The first step is the following 

\blemma
\label{matgpcsa}
\bgloz
  C(\Fz_q [T] / T^n \cdot \Fz_q [T]) \rtimes (\Fz_q [T])^{(n-1)} \rtimes \Fz_q \reg \cong M_{q^n}(\Cz) \otimes C^*(\Fz_q \reg).
\egloz
\elemma

\bproof
Let $e_n \in C(\Fz_q [T] / T^n \cdot \Fz_q [T])$ be the characteristic function of the coset $0 + T^n \cdot \Fz_q [T]$. It is clear that $\menge{v^a e_n v^{-a'}}{a,a' \in (\Fz_q [T])^{(n-1)}}$ are matrix units. Let $e_{a,a'}$ be the canonical rank $1$ operator in $\cL(\ell^2(\Fz_q [T] / T^n \cdot \Fz_q [T]))$ corresponding to the cosets $a + T^n \cdot \Fz_q [T]$ and $a' + T^n \cdot \Fz_q [T]$. We can identify $C(\Fz_q [T] / T^n \cdot \Fz_q [T]) \rtimes (\Fz_q [T])^{(n-1)}$  with $\cL(\ell^2(\Fz_q [T] / T^n \cdot \Fz_q [T])) \cong M_{q^n}(\Cz)$ via 
\bgloz
  v^a e_n v^{-a'} \loma e_{a,a'}. 
\egloz 

Moreover, the action of $\Fz_q \reg$ on $C(\Fz_q [T] / T^n \cdot \Fz_q [T]) \rtimes (\Fz_q [T])^{(n-1)}$ must be inner as we have seen that $C(\Fz_q [T] / T^n \cdot \Fz_q [T]) \rtimes (\Fz_q [T])^{(n-1)}$ is isomorphic to a matrix algebra. The unitaries implementing the action of $\Fz_q \reg$ are given by 
\bgloz
  \sum_{a \in (\Fz_q [T])^{(n-1)}} v^{ba} e_n v^{-a}
\egloz
for $b \in \Fz_q \reg$. 

So on the whole, we obtain the identification 
\bgloz
  C(\Fz_q [T] / T^n \cdot \Fz_q [T]) \rtimes (\Fz_q [T])^{(n-1)} \rtimes \Fz_q \reg \cong \cL(\ell^2(\Fz_q [T] / T^n \cdot \Fz_q [T])) \otimes C^*(\Fz_q \reg)
\egloz
via 
\bgloz
  v^a e_n v^{-a'} t_b \ma e_{a,b^{-1}a'} \otimes V_b
\egloz 
where $V_b$ are the canonical unitary generators of $C^*(\Fz_q \reg)$. 
\eproof

For every character $\chi$ of $\Fz_q \reg$, let $p_{\chi}$ be the spectral projection 
\bgloz
  \tfrac{1}{q-1} \sum_{b \in \Fz_q \reg} \chi(b) t_b 
\egloz 
in $C(\Fz_q [T] / T^n \cdot \Fz_q [T]) \rtimes (\Fz_q [T])^{(n-1)} \rtimes \Fz_q \reg$. As an immediate consequence of Lemma \ref{matgpcsa} we get 

\bcor
\label{K_04}
\bgloz
  K_0(C(\Fz_q [T] / T^n \cdot \Fz_q [T]) \rtimes (\Fz_q [T])^{(n-1)} \rtimes \Fz_q \reg) \cong \bigoplus_{\widehat{\Fz_q \reg}} \Zz \ (\cong \Zz^{q-1}).
\egloz
and free generators for $K_0$ are $\eckl{e_n \cdot p_{\chi}}$, $\chi \in \widehat{\Fz_q \reg}$. 

$K_1(C(\Fz_q [T] / T^n \cdot \Fz_q [T]) \rtimes (\Fz_q [T])^{(n-1)} \rtimes \Fz_q \reg)$ vanishes. 
\ecor
Recall that $e_n \in C(\Fz_q [T] / T^n \cdot \Fz_q [T])$ is the characteristic function of the coset $0 + T^n \cdot \Fz_q [T]$. 

Just a remark on notation: $\eckl{\cdot}$ denotes a class in K-theory. 

By continuity of $K_1$ and with the help of Lemmas \ref{f1}, \ref{f2} and \ref{f3}, we deduce from the previous corollary 
\bcor
\label{K_1}
\bgloz
  K_1(A_{-1}) \cong \gekl{0}.
\egloz
\ecor

It remains to determine $K_0(A_{-1})$. 

\blemma
\label{K_03}
We can identify $K_0(C(\Fz_q [[T]]) \rtimes \Fz_q [T] \rtimes \Fz_q \reg)$ with 
\bgloz
  \Zz[\tfrac{1}{q}] \oplus \bigoplus_{\widehat{\Fz_q \reg} \setminus \gekl{1}} \Zz \ (\cong \Zz[\tfrac{1}{q}] \oplus \Zz^{q-2}). 
\egloz
Moreover, the identification can be chosen so that the $n$-th embedding 
\bgloz
  \iota_n: C(\Fz_q [T] / T^{n} \cdot \Fz_q [T]) \rtimes (\Fz_q [T])^{(n-1)} \rtimes \Fz_q \reg \ri C(\Fz_q [[T]]) \rtimes \Fz_q [T] \rtimes \Fz_q \reg 
\egloz
is given on $K_0$ by 
\bgloz
  (\iota_n)_*(\eckl{e_n}) 
  = 
  \rukl{
  \begin{smallmatrix}
  \tfrac{1}{q^n} \\
  0 \\
  \vdots \\
  0
  \end{smallmatrix}
  }; \ 
  (\iota_n)_*(\eckl{e_n \cdot p_{\chi}}) 
  = 
  \rukl{
  \begin{smallmatrix}
  -\tfrac{1}{q^n} \tfrac{q^n-1}{q-1}\\
  0 \\
  \vdots \\
  1 \\
  \vdots \\
  0
  \end{smallmatrix}
  }  
\egloz
for every $\chi \in \widehat{\Fz_q \reg} \setminus \gekl{1}$. Here, the \an{1} in the image of $\eckl{e_n \cdot p_{\chi}}$ is the entry corresponding to $\chi$ in $\bigoplus_{\widehat{\Fz_q \reg} \setminus \gekl{1}} \Zz$. 

In particular, generators for $K_0(C(\Fz_q [[T]]) \rtimes \Fz_q [T] \rtimes \Fz_q \reg)$ are given by 
\bgloz
  \eckl{\1z_n}, \ n \in \Zz_{\geq 0}
  \text{ and }
  \eckl{p_{\chi}}, \ \chi \in \widehat{\Fz_q \reg} \setminus \gekl{1}.
\egloz
$\1z_n$ is the characteristic function of $T^n \cdot \Fz_q [[T]]$, and $1 \in \widehat{\Fz_q \reg}$ denotes the trivial character.
\elemma

\bproof
With Lemma \ref{f3} in mind, we compute $(\iota_{n,n+1})_*$. By definition,  
\bgloz
  \iota_{n,n+1}(e_n) = \sum_{b \in \Fz_q} v^{bT^n} e_{n+1} v^{-bT^n}; \ \iota_{n,n+1}(t_b) = t_b.
\egloz
Thus, by Corollary \ref{K_04}, we have to determine 
\bgloz
  \eckl{\iota_{n,n+1}(e_n \cdot p_{\chi})} = \eckl{\rukl{\sum_{b \in \Fz_q} v^{bT^n} e_{n+1} v^{-bT^n}} \cdot p_{\chi}}
\egloz
in $K_0(C(\Fz_q [T] / T^{n+1} \cdot \Fz_q [T]) \rtimes (\Fz_q [T])^{(n)} \rtimes \Fz_q \reg)$. 

First of all, we have 
\bgln
\label{MvN}
  && p_{\chi} \cdot \rukl{\sum_{b \in \Fz_q \reg} \psi(b) v^{bT^n} e_{n+1} v^{-bT^n}} \\
  &=& \tfrac{1}{q-1} \sum_{b,b' \in \Fz_q \reg} \chi(b') \psi(b) t_{b'} v^{bT^n} e_{n+1} v^{-bT^n} \nonumber \\
  &=& \tfrac{1}{q-1} \sum_{b,b' \in \Fz_q \reg} \chi(b') \psi(b) v^{b'bT^n} e_{n+1} v^{-b'bT^n} t_{b'} \nonumber \\
  &=& \tfrac{1}{q-1} \sum_{b' \in \Fz_q \reg} \rukl{\sum_{b \in \Fz_q \reg} \psi(b') \psi(b) v^{b'bT^n} e_{n+1} v^{-b'bT^n}} \overline{\psi}(b') \chi(b') t_{b'} 
  \nonumber \\
  &=& \rukl{\sum_{b \in \Fz_q \reg} \psi(b) v^{bT^n} e_{n+1} v^{-bT^n}} \cdot p_{\overline{\psi} \cdot \chi} \nonumber
\egln
for every $\psi$, $\chi$ in $\widehat{\Fz_q \reg}$. This result implies that the projections $(e_n-e_{n+1}) \cdot p_{\overline{\psi} \cdot \chi}$ and $(e_n-e_{n+1}) \cdot p_{\chi} = \rukl{\sum_{b \in \Fz_q \reg} v^{bT^n} e_{n+1} v^{-bT^n}} \cdot p_{\chi}$ are Murray-von Neumann equivalent via the partial isometry 
\bgloz
  p_{\chi} \cdot \rukl{\sum_{b \in \Fz_q \reg} \psi(b) v^{bT^n} e_{n+1} v^{-bT^n}}. 
\egloz
This shows that in $K_0(C(\Fz_q [T] / T^{n+1} \cdot \Fz_q [T]) \rtimes (\Fz_q [T])^{(n)} \rtimes \Fz_q \reg)$, the following equality holds true: 
\bglnoz
  && (q-1) \eckl{e_{n+1}} = \eckl{\sum_{b \in \Fz_q \reg} v^{bT^n} e_{n+1} v^{-bT^n}} = \eckl{e_n-e_{n+1}} \\
  &=& \sum_{\psi \in \widehat{\Fz_q \reg}} \eckl{(e_n-e_{n+1}) \cdot p_{\overline{\psi} \cdot \chi}} = (q-1) \eckl{(e_n-e_{n+1}) \cdot p_{\chi}}
\eglnoz
for every $\chi \in \widehat{\Fz_q \reg}$. Comparing this with 
\bgloz
  (q-1) \eckl{e_{n+1}} = (q-1) \sum_{\psi \in \widehat{\Fz_q \reg}} \eckl{e_{n+1} \cdot p_{\psi}}, 
\egloz
we deduce 
\bgloz
  \eckl{(e_n-e_{n+1}) \cdot p_{\chi}} = \sum_{\psi \in \widehat{\Fz_q \reg}} \eckl{e_{n+1} \cdot p_{\psi}} 
\egloz
for every $\chi \in \widehat{\Fz_q \reg}$. Therefore, 
\bgloz
  (\iota_{n,n+1})_*(\eckl{e_n \cdot p_{\chi}}) = \eckl{\rukl{e_{n+1}+(e_n-e_{n+1})} \cdot p_{\chi}} 
  = \eckl{\rukl{e_{n+1} \cdot p_{\chi}}} + \sum_{\psi \in \widehat{\Fz_q \reg}} \eckl{e_{n+1} \cdot p_{\psi}}.
\egloz
Hence, under the identification 
\bgloz
  K_0(C(\Fz_q [T] / T^{n} \cdot \Fz_q [T]) \rtimes (\Fz_q [T])^{(n-1)} \rtimes \Fz_q \reg) \cong \Zz^{q-1} 
\egloz 
in Corollary \ref{K_04}, we get 
\bgloz
  (\iota_{n,n+1})_* 
  = 
  \rukl{
  \begin{smallmatrix}
  1 & & 0 \\
   & \ddots & \\
  0 & & 1
  \end{smallmatrix}
  }
  +
  \rukl{
  \begin{smallmatrix}
  1 & \dotso & 1 \\
  \vdots & 1 & \vdots \\
  1 & \dotso & 1
  \end{smallmatrix}
  }
  = 
  \rukl{
  \begin{smallmatrix}
  2 & & 1 \\
   & \ddots & \\
  1 & & 2
  \end{smallmatrix}
  }. 
\egloz
Finally, by Lemma \ref{f3}, we compute 
\bgloz
  K_0(C(\Fz_q [[T]]) \rtimes \Fz_q [T] \rtimes \Fz_q \reg) 
  \cong \ilim 
  \gekl{
  \Zz^{q-1}; 
  \rukl{
  \begin{smallmatrix}
  2 & & 1 \\
   & \ddots & \\
  1 & & 2
  \end{smallmatrix}
  }
  }
  \cong \Zz[\tfrac{1}{q}] \oplus \bigoplus_{\widehat{\Fz_q \reg} \setminus \gekl{1}} \Zz.
\egloz
Moreover, we can choose this identification so that $(\iota_n)_*$ is given by 
\bglnoz
  (\iota_n)_*(\eckl{e_n}) 
  = 
  \rukl{
  \begin{smallmatrix}
  \tfrac{1}{q^n} \\
  0 \\
  \vdots \\
  0
  \end{smallmatrix}
  }; \ 
  (\iota_n)_*(\eckl{e_n \cdot p_{\chi}}) 
  = 
  \rukl{
  \begin{smallmatrix}
  -\tfrac{1}{q^n} \tfrac{q^n-1}{q-1}\\
  0 \\
  \vdots \\
  1 \\
  \vdots \\
  0
  \end{smallmatrix}
  }
  \fa \chi \in \widehat{\Fz_q \reg} \setminus \gekl{1}.
\eglnoz

The last statement about the generators of $K_0(C(\Fz_q [[T]]) \rtimes \Fz_q [T] \rtimes \Fz_q \reg)$ follows from the observation that $e_n$ is sent to $\1z_n$ under $\iota_n$. 
\eproof

From now on, we fix this particular description of $K_0(C(\Fz_q [[T]]) \rtimes \Fz_q [T] \rtimes \Fz_q \reg)$.

$\1z_n$ is the characteristic function of $T^n \cdot \Fz_q [[T]] \subseteq \Fz_q [[T]]$. For all $d$ in $\Gamma$,  
\bgloz
  \mu_d(\1z_n) = \chf{T^n \cdot \Fz_q [[T]]}(d^{-1} \sqcup) = \chf{d \cdot (T^n \cdot \Fz_q [[T]])} = \1z_n
\egloz
because every $d$ in $\Gamma$ is invertible in $\Fz_q [[T]]$ ($\mu$ and $\Gamma$ are defined in Section \ref{not}). Moreover, $\mu_d$ certainly leaves $p_{\chi}$ invariant. Therefore, by Lemma \ref{K_03}, $(\mu_d)_* = \id$ on $K_0(C(\Fz_q [[T]]) \rtimes \Fz_q [T] \rtimes \Fz_q \reg)$. This, together with Lemma \ref{f2}, implies 

\bcor
\label{K_02}
\bgloz
  K_0(C(\Fz_q [[T]]) \rtimes (\Gamma \cdot \Fz_q [T]) \rtimes \Fz_q \reg) \cong \Zz[\tfrac{1}{q}] \oplus \bigoplus_{\widehat{\Fz_q \reg} \setminus \gekl{1}} \Zz
\egloz
and the canonical inclusion 
\bgloz
  C(\Fz_q [[T]]) \rtimes \Fz_q [T] \rtimes \Fz_q \reg \ri C(\Fz_q [[T]]) \rtimes (\Gamma \cdot \Fz_q [T]) \rtimes \Fz_q \reg
\egloz
is an isomorphism on $K_0$. 
\ecor

Finally, we have to compute $(\mu_T)_*$ on $K_0(C(\Fz_q [[T]]) \rtimes (\Gamma \cdot \Fz_q [T]) \rtimes \Fz_q \reg)$. We have $\mu_T(\1z_n) = \1z_{n+1}$ and $\mu_T(p_{\chi}) = \1z_1 \cdot p_{\chi}$. Thus, under the identifications in Lemma \ref{K_03} and Corollary \ref{K_02}, we have 
\bgl
\label{mu_T}
  (\mu_T)_*
  = 
  \rukl{
  \begin{smallmatrix}
  \tfrac{1}{q} & -\tfrac{1}{q} & \dotso & -\tfrac{1}{q} \\
  0 & 1 & & 0 \\
  \vdots & & \ddots & \\
  0 & 0 & & 1
  \end{smallmatrix}
  }.
\egl
In particular, $(\mu_T)_*$ is bijective on $K_0(C(\Fz_q [[T]]) \rtimes (\Gamma \cdot \Fz_q [T]) \rtimes \Fz_q \reg)$. Again, combining this result with Lemma \ref{f1} and Corollary \ref{K_02}, we get 

\bcor 
\label{K_01}
\bgloz
  K_0(A_{-1}) \cong \Zz[\tfrac{1}{q}] \oplus \bigoplus_{\widehat{\Fz_q \reg} \setminus \gekl{1}} \Zz
\egloz
and the canonical inclusion 
\bgloz
  C(\Fz_q [[T]]) \rtimes \Fz_q [T] \rtimes \Fz_q \reg \ri A_{-1}
\egloz
is an isomorphism on $K_0$. 

Generators of $K_0(A_{-1})$ are 
\bgloz
  \eckl{\1z_n} \ \hat{=} 
  \rukl{
  \begin{smallmatrix}
  \tfrac{1}{q^n} \\
  0 \\
  \vdots \\
  0
  \end{smallmatrix}
  }
  \text{ and }
  \eckl{\1z \cdot p_{\chi}} \ \hat{=}
  \rukl{
  \begin{smallmatrix}
  0 \\
  \vdots \\
  1 \\
  \vdots \\
  0
  \end{smallmatrix}
  }. 
\egloz
Here, the \an{1} is the entry corresponding to $\chi$ in $\bigoplus_{\widehat{\Fz_q \reg} \setminus \gekl{1}} \Zz$. 
\ecor

So we have obtained a concrete description of the K-theory of 
\bgloz 
  A_{-1} = C_0(\Fz_q ((T))) \rtimes \Fz_q (T) \rtimes \Fz_q \reg. 
\egloz 
We fix this description of $K_0(A_{-1})$ from now on. The next step is to determine the K-theory of 
\bgloz
  A_0 = C_0(\Fz_q ((T))) \rtimes \Fz_q (T) \rtimes (\Fz_q \reg \times \Gamma_0) \cong C_0(\Fz_q ((T))) \rtimes \Fz_q (T) \rtimes \Fz_q \reg \rtimes_{\mu_T} \Zz. 
\egloz

For every $\chi$ in $\widehat{\Fz_q \reg}$, let 
\bgl
\label{x}
  x_{\overline{\chi}} \defeq \sum_{b \in \Fz_q \reg} \overline{\chi}(b) \chf{b+T \cdot \Fz_q [[T]]} \in A_0.
\egl
Moreover, we construct 
\bgl
\label{w}
  w_{\chi} = t_T (\1z \cdot p_1) + x_{\overline{\chi}} p_{\chi} + (\1z \cdot p_{\chi}) t_T^* + (1-\1z \cdot p_1-\1z \cdot p_{\chi})
\egl
in the unitalization $(A_0)^{\sim}$ of $A_0$. $1$ denotes the unit in $(A_0)^{\sim}$. A straightforward computation shows that $w_{\chi}$ is unitary. 

\bprop
\label{KA_0}
We can identify $K_0(A_0)$ with $\bigoplus_{\widehat{\Fz_q \reg} \setminus \gekl{1}} \Zz$ and free generators are $\eckl{\1z \cdot p_{\chi}}$, $1 \neq \chi \in \widehat{\Fz_q \reg}$. 

We also have $K_1(A_0) \cong \bigoplus_{\widehat{\Fz_q \reg} \setminus \gekl{1}} \Zz$ and free generators are $\eckl{w_{\chi}}$, $1 \neq \chi \in \widehat{\Fz_q \reg}$. \eprop

\bproof
$A_0$ can be described as the crossed product $A_{-1} \rtimes_{\mu_T} \Zz$. Thus we can apply the Pimsner-Voiculescu sequence. It looks as follows:
\bgloz
  \gekl{0} \ri K_1(A_0) \overset{\partial}{\ri} K_0(A_{-1}) \overset{\id - (\mu_T)_*}{\lori} K_0(A_{-1}) \ri K_0(A_0) \ri \gekl{0}
\egloz
If we plug in (\ref{mu_T}), then we obtain 
\bglnoz
  && \ker(\id - (\mu_T)_*) = \spkl{\menge{\eckl{\1z \cdot p_1} - \eckl{\1z \cdot p_{\chi}}}{\chi \in \widehat{\Fz_q \reg} \setminus \gekl{1}}}; \\
  && \img(\id - (\mu_T)_*) = \spkl{\menge{\eckl{\1z_n}}{n \in \Zz_{\geq 0}}}.
\eglnoz
As an immediate consequence, we get that $K_0(A_0) \cong \bigoplus_{\widehat{\Fz_q \reg} \setminus \gekl{1}} \Zz$ with free generators $\eckl{\1z \cdot p_{\chi}}$, $1 \neq \chi \in \widehat{\Fz_q \reg}$, as desired. 

To prove our assertion about $K_1$, we have to show that 
\bgloz
  \partial(\eckl{w_{\chi}}) = \eckl{\1z \cdot p_1} - \eckl{\1z \cdot p_{\chi}} \text{ (up to sign)}.
\egloz 
In order to do so, let us have a closer look at the Pimsner-Voiculescu sequence (compare \cite{PV}). It is derived from the Toeplitz extension associated with the crossed product, where the C*-algebra on which $\Zz$ acts is assumed to be unital. As we are in the nonunital case, we have to look at the Toeplitz extension associated to $(A_{-1})^{\sim} \rtimes_{\ti{\mu}_T} \Zz$, i.e. 
\bgl
\label{Toe}
  \gekl{0} \ri \cK \otimes (A_{-1})^{\sim} \ri \cT \ri (A_{-1})^{\sim} \rtimes_{\ti{\mu}_T} \Zz \ri \gekl{0}. 
\egl
Here, $\cK$ is the C*-algebra of compact operators (on some infinite-dimensional separable Hilbert space) and $\cT$ is the C*-subalgebra of $C^*(v) \otimes ((A_{-1})^{\sim} \rtimes_{\ti{\mu}_T} \Zz)$ generated by $v \otimes t_T$ and $\menge{1 \otimes x}{x \in (A_{-1})^{\sim}}$. $C^*(v)$ is the Toeplitz algebra with canonical generator $v$. The quotient map in (\ref{Toe}) maps $v \otimes t_T$ to $t_T$.

Now, to compute $\partial(\eckl{w_{\chi}})$, we consider the partial isometry 
\bgloz
  s_{\chi} 
  = v \otimes t_T (\1z \cdot p_1) + 1 \otimes x_{\overline{\chi}} p_{\chi} + v^* \otimes (\1z \cdot p_{\chi}) t_T^* + 1 \otimes (1-\1z \cdot p_1-\1z \cdot p_{\chi}). 
\egloz
$s_{\chi}$ is mapped to $w_{\chi}$ under the quotient map in (\ref{Toe}). Thus, by definition of $\partial$, we have (up to sign) 
\bglnoz
  && \partial(\eckl{w_{\chi}}) = \eckl{s_{\chi}^* s_{\chi}} - \eckl{s_{\chi} s_{\chi}^*} \\
  &=& \eckl{1 \otimes 1 - (1-vv^*) \otimes (\1z_1 \cdot p_{\chi})} - \eckl{1 \otimes 1 - (1-vv^*) \otimes (\1z_1 \cdot p_1)} \\
  &=& \eckl{(1-vv^*) \otimes (\1z_1 \cdot p_1)} - \eckl{(1-vv^*) \otimes (\1z_1 \cdot p_{\chi})}
\eglnoz
where we used (\ref{MvN}). The last term corresponds to $\eckl{\1z_1 \cdot p_1} - \eckl{\1z_1 \cdot p_{\chi}}$ under the canonical isomorphism $K_0(\cK \otimes (A_{-1})^{\sim}) \cong K_0((A_{-1})^{\sim})$. Finally, using Lemma \ref{K_03}, we deduce that 
\bgloz
  \eckl{\1z_1 \cdot p_1} - \eckl{\1z_1 \cdot p_{\chi}} = (\iota_1)_*(\eckl{p_1}) - (\iota_1)_*(\eckl{p_{\chi}}) 
\egloz
corresponds to 
$
  \rukl{
  \begin{smallmatrix}
  \tfrac{1}{q}(q-1) \\
  -1 \\
  \vdots \\
  -1
  \end{smallmatrix}
  }
  -
  \rukl{
  \begin{smallmatrix}
  -\tfrac{1}{q} \\
  0 \\
  \vdots \\
  1 \\
  \vdots \\
  0
  \end{smallmatrix}
  }
  = 
  \rukl{
  \begin{smallmatrix}
  1 \\
  -1 \\
  \vdots \\
  -1
  \end{smallmatrix}
  }
  -
  \rukl{
  \begin{smallmatrix}
  0 \\
  \vdots \\
  1 \\
  \vdots \\
  0
  \end{smallmatrix}
  }
$
in $\Zz[\tfrac{1}{q}] \oplus \bigoplus_{\widehat{\Fz_q \reg} \setminus \gekl{1}} \Zz$. But by Lemma \ref{K_03}, $\eckl{\1z \cdot p_1} - \eckl{\1z \cdot p_{\chi}}$ also corresponds to 
$
  \rukl{
  \begin{smallmatrix}
  1 \\
  -1 \\
  \vdots \\
  -1
  \end{smallmatrix}
  }
  -
  \rukl{
  \begin{smallmatrix}
  0 \\
  \vdots \\
  1 \\
  \vdots \\
  0
  \end{smallmatrix}
  }
$
in $\Zz[\tfrac{1}{q}] \oplus \bigoplus_{\widehat{\Fz_q \reg} \setminus \gekl{1}} \Zz$, where in the second vector, the \an{1} is the entry corresponding to $\chi$ in $\bigoplus_{\widehat{\Fz_q \reg} \setminus \gekl{1}} \Zz$. 

Thus, $\eckl{\1z_1 \cdot p_1} - \eckl{\1z_1 \cdot p_{\chi}} = \eckl{\1z \cdot p_1} - \eckl{\1z \cdot p_{\chi}}$ in $K_0(A_{-1})$. This proves our claim.
\eproof

At this point, we remark that $A_0$ can be described as a Cuntz-Krieger algebra. This leads to an alternative way of computing the K-theory for $A_0$. 

First of all, $C(\Fz_q [[T]]) \rtimes \Fz_q [T] \rte_{\mu_T} \Nz$ is generated by the isometries $v^a t_T$, $a \in \Fz_q$, whose range projections sum up to $1$. Here $t_T$ is the isometry which implements the endomorphism $\mu_T$. Thus $C(\Fz_q [[T]]) \rtimes \Fz_q [T] \rte_{\mu_T} \Nz$ is isomorphic to the Cuntz algebra $\cO_q$ (by the universal property of $\cO_q$ and since the Cuntz algebra is simple, see \cite{Cun1}). 

Secondly, consider the crossed product $\cO_q \rtimes (\Zz / (q-1) \Zz)$ with respect to the action 
\bgloz
  S_j 
  \loma
  \bfa
    S_j \fuer j=1 \\
    \zeta^{j-2} S_j \falls j \geq 2
  \efa
\egloz
for a primitive $(q-1)$-th root of unity $\zeta$, where the $S_j$ are the canonical generators of $\cO_q$. We claim that 
\bgloz
  \cO_q \rtimes (\Zz / (q-1) \Zz) \cong C(\Fz_q [[T]]) \rtimes \Fz_q [T] \rte_{\mu_T} \Nz \rtimes \Fz_q \reg.
\egloz
To show this, choose a generator $b$ of $\Fz_q \reg$ and consider the isometries 
\bgloz
  t_T \text{ and } \tfrac{1}{\sqrt{q-1}} \sum_{n=0}^{q-2} (\zeta^{2-j})^n v^{(b^n)} t_T \fuer 2 \leq j \leq q.
\egloz
These $q$ isometries generate $C(\Fz_q [[T]]) \rtimes \Fz_q [T] \rte_{\mu_T} \Nz$, and their range projections sum up to $1$. Moreover, we have $t_b t_T t_b^* = t_T$ and 
\bgloz
  t_b \rukl{\tfrac{1}{\sqrt{q-1}} \sum_{n=0}^{q-2} (\zeta^{2-j})^n v^{(b^n)} t_T} t_b^* 
  = \zeta^{j-2} \rukl{\tfrac{1}{\sqrt{q-1}} \sum_{n=0}^{q-2} (\zeta^{2-j})^n v^{(b^n)} t_T}
\egloz
for $2 \leq j \leq q$. Thus we have found a $\Zz / (q-1) \Zz \cong \Fz_q \reg$-invariant isomorphism 
\bgloz
  \cO_q \cong C(\Fz_q [[T]]) \rtimes \Fz_q [T] \rte_{\mu_T} \Nz.
\egloz
Here we again used the universal property of $\cO_q$ together with the fact that $\cO_q$ is simple (see \cite{Cun1}). We conclude that 
\bgloz
  \cO_q \rtimes (\Zz / (q-1) \Zz) \cong C(\Fz_q [[T]]) \rtimes \Fz_q [T] \rte_{\mu_T} \Nz \rtimes \Fz_q \reg, 
\egloz
as claimed. 

And thirdly, by \cite{CuEv} we know that $\cO_q \rtimes (\Zz / (q-1) \Zz)$ is isomorphic to the Cuntz-Krieger algebra $\cO_A$ associated with the matrix
$
  A
  = 
  \rukl{
  \begin{smallmatrix}
  2 & & 1 \\
  & \ddots & \\
  1 & & 2 
  \end{smallmatrix}
  }
$. 

Thus we get 
\bgloz
  C(\Fz_q [[T]]) \rtimes \Fz_q [T] \rte_{\mu_T} \Nz \rtimes \Fz_q \reg \cong \cO_A.
\egloz
Now this isomorphism can be worked out explicitly, and we can compute the K-theory of $\cO_A$ in an explicit way (compare \cite{CuKr} and \cite{Cun2}). Therefore we obtain a concrete description for the K-theory of $C(\Fz_q [[T]]) \rtimes \Fz_q [T] \rte_{\mu_T} \Nz \rtimes \Fz_q \reg$. 

Finally, it follows from our computations that the canonical homomorphism 
\bgloz
  C(\Fz_q [[T]]) \rtimes \Fz_q [T] \rte_{\mu_T} \Nz \rtimes \Fz_q \reg \ri A_0
\egloz
is an isomorphism on K-theory. So this is an alternative route of computing the K-theory of $A_0$.

\section{Commuting unitaries}
\label{communit}

Now we come to the crucial point in our computations. We have invested some effort in describing the K-groups of $A_0$ as explicitly as possible. The reason is that we are interested in the following observation: 
\blemma
\label{comm}
For every $i$ and $1 \neq \chi \in \widehat{\Fz_q \reg}$, we have 
\bgloz
  \mu_{f_i}(\1z \cdot p_{\chi}) = \1z \cdot p_{\chi} \text{ and } \mu_{f_i}(w_{\chi}) = w_{\chi}.
\egloz
\elemma
Recall that $f_1, f_2, f_3, \dotsc$ is an enumeration of the irreducible polynomials in $\Fz_q[T]$ with constant term 1.

\bproof
We have $\mu_{f_i}(\1z) = \chf{f_i \cdot \Fz_q [[T]]} = \1z$ as $f_i$ is invertible in $\Fz_q [[T]]$. This shows $\mu_{f_i}(\1z \cdot p_{\chi}) = \1z \cdot p_{\chi}$. 

To show that $w_{\chi}$ (defined in (\ref{w})) is $\mu_{f_i}$-invariant, it remains to prove that $x_{\overline{\chi}}$ is $\mu_{f_i}$-invariant. By construction of $x_{\overline{\chi}}$ (defined in (\ref{x})), it suffices to prove that $\chf{b+T \cdot \Fz_q [[T]]}$ is $\mu_{f_i}$-invariant for all $b$ in $\Fz_q \reg$. Since $\mu_{f_i}(\chf{b+T \cdot \Fz_q [[T]]}) = \chf{f_i \cdot (b+T \cdot \Fz_q [[T]])}$, we have to show that $b+T \cdot \Fz_q [[T]] = f_i \cdot (b+T \cdot \Fz_q [[T]])$. As $f_i$ has constant term $1$, it is clear that \an{$\subseteq$} holds. To prove the reverse inclusion, take an arbitrary element $b+Tx$ in $b+T \cdot \Fz_q [[T]]$. Then 
\bgloz
  b+Tx = f_i \cdot (b + \underbrace{f_i^{-1}}_{\in \Fz_q [[T]]} \cdot (\underbrace{(1-f_i)b+Tx}_{\in T \cdot \Fz_q [[T]]})) \in f_i \cdot (b+T \cdot \Fz_q [[T]]).
\egloz
This proves our lemma.
\eproof

As we will see, this simple observation plays a very important role in our computations. Moreover, note that this observation heavily relies on the fact that we have applied our duality theorem to pass from the finite adele ring to the infinite one. The reason why our lemma holds true basically is that all the $f_i$ are invertible in $\Fz_q [[T]]$. But this only happens in the canonical subring of the infinite adele ring, whereas in the canonical subring of the finite adele ring, there is for each polynomial $f_i$ a finite place where $f_i$ is not invertible. 

Now, the reason why this observation is so important is that it allows us to produce generators for the K-theory of $A_m$. 

We fix the following notations: Let $t(i)$ be the unitary $t_{f_i}$ in the multiplier algebra of $A_m$ for every $1 \leq i \leq m$. We know that $A_m \cong A_{m-1} \rtimes_{\mu_{f_m}} \Zz$, and we denote by $\partial_m$ the boundary map in the corresponding Pimsner-Voiculescu sequence. It will become clear from the context whether we mean the index map or the exponential map. Let 
\bgl
\label{1t}
  (\1z \cdot p_{\chi}, t(m)) \defeq t(m) (\1z \cdot p_{\chi}) + (1-\1z \cdot p_{\chi}) \in (A_m)^{\sim}.
\egl 
Here $1$ is the unit in $(A_m)^{\sim}$. 

\blemma
$(\1z \cdot p_{\chi}, t(m))$ is a unitary in $(A_m)^{\sim}$ with 
\bgl
\label{d_m0}
  \partial_m(\eckl{(\1z \cdot p_{\chi}, t(m))}) = \eckl{\1z \cdot p_{\chi}} \in K_0(A_{m-1})
\egl
(up to sign). 
\elemma

\bproof
First of all, $(\1z \cdot p_{\chi}, t(m))$ is a unitary since $t(m)$ commutes with $\1z \cdot p_{\chi}$ as $\1z \cdot p_{\chi}$ is $\mu_{f_m}$-invariant. To prove (\ref{d_m0}), we have to look at the Toeplitz extension (with generalized Toeplitz algebra $\cT$) associated to $(A_{m-1})^{\sim} \rtimes_{\ti{\mu}_{f_m}} \Zz$ as in the proof of Proposition \ref{KA_0}. Here $\ti{\mu}_{f_m}$ is the extension of $\mu_{f_m}$ to the unitalization. 

We find that the partial isometry $\ti{s}_{\chi} \defeq v \otimes t(m) (\1z \cdot p_{\chi}) + 1 \otimes (1-\1z \cdot p_{\chi})$ is mapped to $(\1z \cdot p_{\chi}, t(m))$ under the quotient map $\cT \ri (A_{m-1})^{\sim} \rtimes_{\ti{\mu}_{f_m}} \Zz$. Thus, 
\bgloz
  \partial_m(\eckl{(\1z \cdot p_{\chi}, t(m))}) = \eckl{\ti{s}_{\chi}^* \ti{s}_{\chi}} - \eckl{\ti{s}_{\chi} \ti{s}_{\chi}^*} = \eckl{(1-vv^*) \otimes \1z \cdot p_{\chi}}
\egloz
(up to sign) and the last term corresponds to $\eckl{\1z \cdot p_{\chi}}$ under the canonical identification $K_0(\cK \otimes (A_{m-1})^{\sim}) \cong K_0((A_{m-1})^{\sim})$. This proves our lemma. 
\eproof

In the following, we produce generators for $K_*(A_m)$ by comparing our situation with higher-dimensional commutative tori. We denote by $K_*(A_m)$ the $\Zz / 2 \Zz$-graded abelian group $K_0(A_m) \oplus K_1(A_m)$. 

For each $l \in \Zz \pos$, let $z_0, \dotsc, z_l$ be the canonical unitary generators of $C(\Tz^{l+1})$. Choose some $1 \neq \chi \in \widehat{\Fz_q \reg}$. Let $\Gamma'_m$ be the subgroup of $\Gamma$ generated by the polynomials $f_1, \dotsc, f_m$, i.e. 
\bgl
\label{Gamma'}
  \Gamma'_m \defeq \spkl{f_1, \dotsc, f_m}.
\egl
By universal property of $C(\Tz^{l+1})$, the commuting unitaries $w_{\chi}, t(i_1), \dotsc, t(i_l)$ (for some $1 \leq i_1 < \dotsb < i_l \leq m$) give rise to a homomorphism 
\bgloz
  C(\Tz^{l+1}) \ri (A_0)^{\sim} \rtimes_{\ti{\mu}} \Gamma'_m; \ z_0 \ma w_{\chi}, z_j \ma t(i_j). 
\egloz
Here $\ti{\mu}$ is the extension of $\mu$ to the unitalization. Note that we can construct such a homomorphism precisely because of Lemma \ref{comm}. 

We denote by $\eckl{w_{\chi}, t(i_1), \dotsc, t(i_l)}$ the image of $\eckl{z_0} \times \dotsb \times \eckl{z_l}$ (see \cite{HiRo}, 4.7 for the definition of the product on K-theory) under this homomorphism in K-theory. A priori, $\eckl{w_{\chi}, t(i_1), \dotsc, t(i_l)}$ lies in $K_*((A_0)^{\sim} \rtimes_{\ti{\mu}} \Gamma'_m)$. However, we observe the following: 

\blemma
\label{d_m1}
$\eckl{w_{\chi}, t(i_1), \dotsc, t(i_l)}$ lies in 
\bgloz
  K_*(A_m) = K_*(A_0 \rtimes_{\mu} \Gamma'_m) \cong \ker \rukl{K_*((A_0)^{\sim} \rtimes_{\ti{\mu}} \Gamma'_m) \ri K_*(C^*(\Gamma'_m))}.
\egloz
\elemma

\bproof
The identification 
\bgloz
  K_*(A_0 \rtimes_{\mu} \Gamma'_m) \cong \ker \rukl{K_*((A_0)^{\sim} \rtimes_{\ti{\mu}} \Gamma'_m) \ri K_*(C^*(\Gamma'_m))}
\egloz
is justified by the split-exact sequence 
\bgloz
  \gekl{0} \ri A_0 \rtimes_{\mu} \Gamma'_m \ri (A_0)^{\sim} \rtimes_{\ti{\mu}} \Gamma'_m \ri C^*(\Gamma'_m) \ri \gekl{0}.
\egloz

Now, consider the commutative diagram
\bgloz
  \begin{CD}
  C_0(\Rz) \otimes C(\Tz^l) @>>> C(\Tz^{l+1}) @>>> C(\Tz^l) \\
  @. @VVV @VVV \\
  @. (A_0)^{\sim} \rtimes_{\ti{\mu}} \Gamma'_m @>>> C^*(\Gamma'_m)
  \end{CD}
\egloz
The first row is split-exact. Moreover, $\eckl{z_0} \times \dotsb \times \eckl{z_l}$ clearly comes from $K_*(C_0(\Rz) \otimes C(\Tz^l))$. Thus, $\eckl{z_0} \times \dotsb \times \eckl{z_l}$ is mapped to $0$ under the homomorphism $K_*(C(\Tz^{l+1})) \ri K_*(C(\Tz^l))$. Since the diagram above commutes, $\eckl{w_{\chi}, t(i_1), \dotsc, t(i_l)}$ must lie in the kernel of $K_*((A_0)^{\sim} \rtimes_{\ti{\mu}} \Gamma'_m) \ri K_*(C^*(\Gamma'_m))$.
\eproof

\blemma
If $i_l = m$, $l \geq 1$, then 
\bgl
\label{d_m2}
  \partial_m(\eckl{w_{\chi}, t(i_1), \dotsc, t(i_l)}) = \eckl{w_{\chi}, t(i_1), \dotsc, t(i_{l-1})} \text{ (up to sign)}. 
\egl

\elemma

\bproof
Under the boundary map $K_*(C(\Tz^{l+1})) \ri K_{*+1}(C(\Tz^l))$ associated with the Toeplitz extension of $C(\Tz^{l+1}) \cong C(\Tz^l) \rtimes_{\id} \Zz$, $\eckl{z_0} \times \dotsb \times \eckl{z_l}$ is mapped to $\eckl{z_0} \times \dotsb \times \eckl{z_{l-1}}$ (up to sign). Therefore, by naturality of the Pimsner-Voiculescu sequence, our claim follows.
\eproof

Similarly, we can consider $C(\Tz^l)$ with canonical unitary generators $\ti{z}_i$, $1 \leq i \leq l$ and the homomorphism $C(\Tz^l) \ri (A_0)^{\sim} \rtimes_{\ti{\mu}} \Gamma'_m$ ($\Gamma'_m$ is defined in (\ref{Gamma'})) given by 
\bgloz
  \ti{z}_1 \ma (\1z \cdot p_{\chi}, t(i_1)); \ \ti{z}_j \ma t(i_j) \fuer 2 \leq j \leq l, 1 \leq i_1 < \dotsb < i_l \leq m. 
\egloz
The unitary $(\1z \cdot p_{\chi}, t(i_1))$ is defined as in (\ref{1t}). Again, this homomorphism exists because the unitaries $(\1z \cdot p_{\chi}, t(i_1)), t(i_2), \dotsc, t(i_l)$ commute (see Lemma \ref{comm}). 

Let $\eckl{\1z \cdot p_{\chi}, t(i_1), t(i_2), \dotsc, t(i_l)}$ be the image of $\eckl{\ti{z}_1} \times \dotsb \times \eckl{\ti{z}_l}$ under the homomorphism above in K-theory. In complete analogy to the preceding two lemmas, we get 
\blemma
\label{d_m3}
$\eckl{\1z \cdot p_{\chi}, t(i_1), t(i_2), \dotsc, t(i_l)}$ lies in 
\bgloz
  K_*(A_m) = K_*(A_0 \rtimes_{\mu} \Gamma'_m) \cong \ker \rukl{K_*((A_0)^{\sim} \rtimes_{\ti{\mu}} \Gamma'_m) \ri K_*(C^*(\Gamma'_m))} 
\egloz
\elemma
and 
\blemma
\bgl
\label{d_m4}
  \partial_m(\eckl{\1z \cdot p_{\chi}, t(i_1), t(i_2), \dotsc, t(i_l)}) = \eckl{\1z \cdot p_{\chi}, t(i_1), t(i_2), \dotsc, t(i_{l-1})}
\egl
(up to sign) for $i_l=m$ and $l>1$ .
\elemma

\section{Final result}
\label{final}

Now we are ready to compute the K-theory of $A_m$ (see Section \ref{not} for the definition of $A_m$). We just have to put everything together.

\bprop
\label{KA_m}
We have (with $\Gamma'_m$ defined in (\ref{Gamma'})) 
\bgloz
  K_*(A_m) \cong K_*(A_0) \otimes_{\Zz} \extalg(\Gamma'_m). 
\egloz

Generators for $K_0$ are $\eckl{w_{\chi}, t(i_1), \dotsc, t(i_k)}$ for $1 \leq i_1 < \dotsb < i_k \leq m$, $1 \leq k$ odd; 
$\eckl{\1z \cdot p_{\chi}, t(i_1), \dotsc, t(i_l)}$ for $1 \leq i_1 < \dotsb < i_l \leq m$, $0 \leq l$ even; with $\chi \in \widehat{\Fz_q \reg} \setminus \gekl{1}$. 

Generators for $K_1$ are $\eckl{w_{\chi}, t(i_1), \dotsc, t(i_k)}$ for $1 \leq i_1 < \dotsb < i_k \leq m$, $0 \leq k$ even; 
$\eckl{\1z \cdot p_{\chi}, t(i_1), \dotsc, t(i_l)}$ for $1 \leq i_1 < \dotsb < i_l \leq m$, $1 \leq l$ odd; with $\chi \in \widehat{\Fz_q \reg} \setminus \gekl{1}$. 
\eprop

Recall that $K_*(A_m)$ is the $\Zz / 2 \Zz$-graded abelian group $K_0(A_m) \oplus K_1(A_m)$. The isomorphism in this proposition is meant as an isomorphism of $\Zz / 2 \Zz$-graded abelian groups, where $\extalg(\Gamma'_m)$ is canonically graded and we consider graded tensor products. 

\bproof
First of all, the statement makes sense because of Lemma \ref{d_m1} and Lemma \ref{d_m3}. Now, to prove our claim, we proceed inductively. For $m=0$ the claim about the generators has been proven in Proposition \ref{KA_0}. 

Assume that $m \geq 1$ and that our assertion holds true for $m-1$. We consider the Pimsner-Voiculescu sequence associated to $A_m \cong A_{m-1} \rtimes_{\mu_{f_m}} \Zz$. The boundary map $\partial_m$ is surjective because we have (up to sign) 
\bgloz
  \partial_m(\eckl{w_{\chi}, t(i_1), \dotsc, t(i_k), t(m)}) = \eckl{w_{\chi}, t(i_1), \dotsc, t(i_k)}
\egloz
for every $i_1 < \dotsc < i_k < m, 0 \leq k$ by (\ref{d_m2}) and 
\bgloz
  \partial_m(\eckl{\1z \cdot p_{\chi}, t(i_1), \dotsc, t(i_l), t(m)}) = \eckl{\1z \cdot p_{\chi}, t(i_1), \dotsc, t(i_l)}
\egloz
for every $i_1 < \dotsc < i_l < m, 0 \leq l$ by (\ref{d_m0}) and (\ref{d_m4}).

Therefore, our claim follows from the exactness of the Pimsner-Voiculescu sequence for $A_m \cong A_{m-1} \rtimes_{\mu_{f_m}} \Zz$ and by the induction hypothesis. In particular, we have $(\mu_{f_m})_* = \id$ on $K_*(A_{m-1})$ for every $m \in \Zz \pos$.
\eproof

Finally, this result allows us to compute the K-theory of the ring C*-algebra associated to $\Fz_q[T]$. By our considerations in Section \ref{dual}, we know that the K-theory of the ring C*-algebra $\fA [\Fz_q[T]]$ can be identified with the K-theory of $C_0(\Fz_q ((T))) \rtimes \Fz_q (T) \rtimes \Fz_q (T) \reg$. Moreover, we have
\bgloz
  C_0(\Fz_q ((T))) \rtimes \Fz_q (T) \rtimes \Fz_q (T) \reg \cong \ilim A_m. 
\egloz
Thus, using continuity of K-theory together with Propositions \ref{KA_0} and \ref{KA_m}, we arrive at the following final result: 

\btheo
\label{Kff}
$K_*(\fA [\Fz_q[T]]) \cong \ti{K}_0(C^*(\Fz_q \reg)) \otimes_{\Zz} \extalg(\Gamma)$. 
\etheo

$\ti{K}_0(C^*(\Fz_q \reg))$ denotes the reduced K-theory of $C^*(\Fz_q \reg)$, i.e. the cokernel of the canonical map $K_0(\Cz) \ri K_0(C^*(\Fz_q \reg))$. $\Gamma$ is defined in Section \ref{not}. Moreover, $K_*(\fA [\Fz_q[T]])$ is the $\Zz / 2 \Zz$-graded abelian group $K_0(\fA [\Fz_q[T]]) \oplus K_1(\fA [\Fz_q[T]])$, and the isomorphism in the theorem above is meant as an isomorphism of $\Zz / 2 \Zz$-graded abelian groups. Here $\ti{K}_0(C^*(\Fz_q \reg))$ is trivially graded, $\extalg(\Gamma)$ is canonically graded and we consider graded tensor products.

Our computations show how to define a product structure on $K_*(\fA [\Fz_q[T]])$ which corresponds to the canonical product structure on  $\ti{K}_0(C^*(\Fz_q \reg)) \otimes_{\Zz} \extalg(\Gamma)$ under the isomorphism above. Actually, it follows from Lemma \ref{comm} that for every $\chi$ in $\widehat{\Fz_q \reg} \setminus \gekl{1}$, the elements $t(1)(\1z \cdot p_{\chi}), t(2)(\1z \cdot p_{\chi}), t(3)(\1z \cdot p_{\chi}), \dotsc$ are commuting unitaries in $(\1z \cdot p_{\chi})(C_0(\Fz_q ((T))) \rtimes \Fz_q (T) \rtimes \Fz_q (T) \reg)(\1z \cdot p_{\chi})$. So they give rise to a homomorphism of the algebra of continuous functions on the infinite dimensional torus to $C_0(\Fz_q ((T))) \rtimes \Fz_q (T) \rtimes \Fz_q (T) \reg$. It follows from Proposition \ref{KA_m} that this homomorphism induces an embedding on K-theory. Thus we just have to carry over the product structure on the K-theory of the infinite dimensional torus to $K_*(\fA [\Fz_q[T]]) \cong K_*(C_0(\Fz_q ((T))) \rtimes \Fz_q (T) \rtimes \Fz_q (T))$. 

As the last comment, we point out that there are striking similarities between the number field case and the function field case (compare Theorem \ref{Knf} and Theorem \ref{Kff}). So from this point of view, our results fit nicely into the general picture concerning analogies between number fields and function fields.


\begin{thebibliography}{99}

\bibitem[Cun1]{Cun1} J. \textsc{Cuntz},
	\emph{Simple C*-algebras generated by isometries},
	Comm. Math. Phys. \emph{57} (1977), 173-185. 

\bibitem[Cun2]{Cun2} J. \textsc{Cuntz},
	\emph{A class of C*-algebras and topological Markov chains II: Reducible Chains and the Ext-functor for C*-algebras},
  Invention. Math. \emph{63} (1981), 25-40.

\bibitem[Cun3]{Cun3} J. \textsc{Cuntz},
	\emph{C*-algebras associated with the $ax + b$-semigroup over $\mathbb{N}$} (English), 
	Corti\~{n}as, Guillermo (ed.) et al., K-theory and noncommutative geometry. 
	Proceedings of the ICM 2006 satellite conference, Valladolid, Spain, August 31-September 6, 2006. 
	Z\"{u}rich: European Mathematical Society (EMS). Series of Congress Reports, 201-215 (2008).

\bibitem[CuEv]{CuEv} J. \textsc{Cuntz} and D. E. \textsc{Evans}, 
  \emph{Some remarks on the C*-algebras associated with certain topological Markov chains}, 
  Math. Scand. \emph{48} (1981), 235-240.

\bibitem[CuKr]{CuKr} J. \textsc{Cuntz} and W. \textsc{Krieger}, 
  \emph{A class of C*-algebras and topological Markov chains}, 
  Invention. Math. \emph{56} (1980), 251-268.

\bibitem[CuLi1]{CuLi1} J. \textsc{Cuntz} and X. \textsc{Li}, 
  \emph{The Regular C*-algebra of an Integral Domain}, arXiv: 0807.1407, 
  to appear in the proceedings of the conference in honour of A. Connes' 60th birthday.

\bibitem[CuLi2]{CuLi2} J. \textsc{Cuntz} and X. \textsc{Li}, 
  \emph{C*-algebras associated with integral domains and crossed products by actions on adele spaces}, 
  arXiv:0906.4903, 
  to appear in the Journal of Noncommutative Geometry.

\bibitem[HiRo]{HiRo} N. \textsc{Higson} and J. \textsc{Roe},
	\emph{Analytic K-Homology}, Oxford University Press, New York, 2000.

\bibitem[Li]{Li} X. \textsc{Li},
	\emph{Ring C*-algebras}, arXiv:0905.4861, preprint.

\bibitem[PV]{PV} M. \textsc{Pimsner} and D. \textsc{Voiculescu},
	\emph{Exact sequences for K-groups and Ext-groups of certain cross-product C*-algebras}, J. Operator Theory \emph{4} (1980), 93-118.
				
\end{thebibliography}
\end{document}